\documentclass[12pt]{amsart}
\usepackage{graphicx}
\allowdisplaybreaks[2]
\newtheorem{theorem}{Theorem}[section]
\newtheorem{corollary}[theorem]{Corollary}

\newtheorem{proposition}[theorem]{Proposition}
\theoremstyle{definition}

\theoremstyle{remark}

\numberwithin{equation}{section}

\newcommand{\norm}[1]{\Vert #1 \Vert}

\newcommand{\B}{\mathcal{B}}

\newcommand{\U}{\mathcal{U}}
\newcommand{\C}{\mathbb{C}}
\newcommand{\N}{\mathbb{N}}
\newcommand{\R}{\mathbb{R}}
\newcommand{\borel}{\Sigma_{\mathbb{C}}}
\newcommand{\one}{\mathbf{1}}
\newcommand{\lam}{\boldsymbol{\lambda}}
\newcommand{\take}{\!\setminus\!}
\newcommand{\<}{\left<}
\renewcommand{\>}{\right>}

\newcommand{\AC}{{\rm AC}}
\newcommand{\BV}{{\rm BV}}
\newcommand{\supp}{{\rm supp}}
\newcommand{\trho}{\tilde{\rho}}

\begin{document}

\title{Functional calculus extensions on dual spaces}
\author{Venta Terauds}
\address{School of Mathematical and Physical Sciences\\
University of Newcastle\\
University Drive\\
Callaghan NSW 2308 Australia} \email{Venta.Terauds@newcastle.edu.au}

\thanks{The author would like to acknowledge the support of an Australian Postgraduate Award and to thank Ian Doust for his comments on an early draft of this paper.}


%\keywords{}%

%\date{\today}
%\dedicatory{}%
%\commby{}%
% ----------------------------------------------------------------

\begin{abstract}
In this note, we show that if a Banach space $X$ has a predual, then every bounded linear operator on $X$ with a continuous functional 
calculus admits a bounded Borel functional calculus. A consequence of this is that on such a Banach space, the classes of finitely spectral
and prespectral operators coincide. We also apply this result to give some sufficient conditions for an
operator with an absolutely continuous functional calculus to admit a bounded Borel one.

\end{abstract}

\maketitle

\section{Introduction}

If it is known that a bounded linear operator $T$ on a Banach space $X$ has a $\U$ functional calculus for some 
algebra of functions $\U$, it is a natural question to ask if there exists a larger algebra of functions, 
some $\mathcal{W}$ containing $\U$, such that $T$ has a $\mathcal{W}$ functional calculus. 
By saying that the operator $T\in L(X)$ has a $\U$ functional calculus for a Banach algebra
of functions $\U$, we mean that there exists a bounded algebra homomorphism $\Psi: \U \to L(X)$ such that $\Psi(\one) = I$ and $\Psi(\lam) = T$.
Implicit here is that $\U = \U(\sigma)$, where $\sigma$ is a (non-empty) compact subset of the plane, and that $\U$ contains the polynomials defined on $\sigma$. The functions $\one$ and $\lam$ refer to the constant function one and the identity function on $\sigma$. As long as $\U$ contains the rational functions with 
poles off $\sigma$, then whenever $T$ has a $\U$ functional calculus it is immediate that the spectrum of $T$ is contained in $\sigma$. 
If an algebra contains sufficiently many idempotent functions, then its functional calculus will provide a family of invariant
subspaces for $T$, and often an integral representation of $T$ with respect to a family of projections.

The algebra $\B(\sigma)$, consisting of all bounded Borel-measurable functions defined on a set $\sigma\subseteq\C$, is one such ``idempotent-rich'' algebra. When an operator admits a $C(\sigma)$ functional calculus, where $C(\sigma)$ is the algebra of continuous functions defined on a set $\sigma\subseteq\C$, one thus looks to extend this to $\B(\sigma)$. Often the structure of the space $X$ ensures the availability of such an extension. For example, if $X$ is a reflexive Banach space, or more generally, one containing no copy of $c_0$, then each operator $T\in L(X)$ with a $C(\sigma)$ functional calculus admits a $\B(\sigma(T))$ functional calculus.   
The latter result may be found in \cite{DdeLaub}, along with several others that specify classes of Banach spaces on which such extensions of a $C(\sigma)$ functional calculus are (and are not) assured. 

The authors note in \cite{DdeLaub} that ``the situation when 
$X = \ell^\infty$ seems to be open''. In Section \ref{section_context}, we answer this question, using a technique due to Sills \cite{sills} 
to show that if $X$ is a 
dual Banach space, then a $C(\sigma)$ functional calculus for an operator $T \in L(X)$ 
may always be extended to a $\B(\sigma)$ functional calculus for $T$. As is usual in the study of such operators, we frame this result in the context of the rich theory of spectral, prespectral and finitely spectral operators. 

A different application of this functional calculus extension is given in the third section of the paper, where we consider real $\AC(\sigma)$ operators. The algebra $\AC(\sigma)$ contains the absolutely continuous functions defined on a non-empty compact set $\sigma\subseteq\C$, and 
an $\AC(\sigma)$ operator is one that possesses an $\AC(\sigma)$ functional calculus. When $\sigma\subseteq\R$ we give sufficient conditions for an $\AC(\sigma)$ operator to possess a $\B(\sigma)$ functional calculus.

\section{$C(\sigma)$ functional calculus}\label{section_context}

A scalar-type spectral operator is one that may be represented as an integral with respect to a spectral measure, that is, a projection valued measure that is countably additive in the strong operator topology. This measure is called a resolution of the identity for the operator. (For full definitions of spectral measures and spectral operators, we refer the reader to \cite{dunford3} or \cite{dowson}.) By weakening the additivity requirement
on the measure, one gains the classes of prespectral and finitely spectral operators (and their scalar types): prespectral operators
admit a resolution of the identity that is countably additive in a weaker operator topology \cite{dunfordspec} and finitely spectral operators
a finitely additive resolution of the identity \cite{nagy}. The scalar-type operators of each class are those with an integral representation. 

It is a standard technique of spectral theory to construct an integral representation for an operator and then use this to define an extended functional calculus. For example if $X$ is a Banach space containing no copy of $c_0$, every map from $C(\sigma)$ to $X$ is weakly compact. This property ensures that a map defining a $C(\sigma)$ functional calculus for an operator $T\in L(X)$ can be used to construct a spectral measure and an integral representation for $T$ \cite{spain}. The operator $T$ is thus scalar-type spectral, and admits a $\B(\sigma(T))$ functional calculus that is defined via its integral representation.

We proceed here in the opposite direction. 
Given an operator $T$ on a dual Banach space $X$ with a $C(\sigma)$ functional calculus, we directly extend this to $\B(\sigma)$, and then use
the extended functional calculus to construct a spectral measure which will provide an integral representation for the operator. The extension method is due to Sills, who proved in \cite[Theorem 6.1]{sills} that a continuous homomorphism from a commutative Banach algebra $\U$ to $L(X)$ may be extended to a homomorphism from $\U^{**}$ to $L(X)$ when $X$ is reflexive. He further stated the following theorem \cite[Theorem 6.2]{sills}, which is proved by essentially the same technique; we provide details of the proof here as we shall use the construction in what follows.
Firstly, recall the following. Given an algebra of functions $\U$, the second dual $\U^{**}$ is itself an algebra under Arens multiplication, which we shall denote by $\star$ and define as in \cite{sills}. For fixed $\varphi \in \U^{**}$, $\psi \mapsto \psi \star \varphi$ is weak-$\ast$ continuous on $\U^{**}$, and for fixed $f\in \U$ (considering $f$ as an element of $\U^{**}$), $\psi \mapsto f \star \psi$ is weak-$\ast$ continuous  on $\U^{**}$.

\begin{theorem}\label{sillsextn}
Let $X$ be a Banach space such that $X = Y^*$ for some Banach space $Y$. Let $\mathcal{U}$ be a commutative Banach algebra and 
$\rho_0: \mathcal{U} \to L(X)$ be a norm-continuous algebra homomorphism. Then $\rho_0$ can be extended to a homomorphism 
$\rho : \mathcal{U}^{**} \to L(X)$ such that $\norm{\rho} = \norm{\rho_0}$ and $\rho$ is continuous from $\mathcal{U}^{**}$ with the 
weak-$\ast$ topology into $L(X)$ with the weak-$\ast$ operator topology. 
%The homomorphism $\rho$ is the unique such extension of $\rho_0$ with this continuity property.
\end{theorem}
\begin{proof}
For $y\in Y, x\in X$, the functional on $\U$ defined by $f\mapsto \< y, \rho_0(f)x\>$ is linear and it is bounded by $\norm{y}\norm{x}\norm{\rho_0}$. Thus there exists $\gamma_{y,x} \in \U^*$ such that $\<f,\gamma_{y,x}\> = \< y, \rho_0(f)x\>$ for all $f\in\U$. 
Then for each pair $\varphi\in\mathcal{U}^{**}$, $x\in X$, the functional on $Y$ defined by $y\mapsto \<\gamma_{y,x}, \varphi\>$ is linear and is bounded by $\norm{x}\norm{\rho_0}\norm{\varphi}$. Hence there exists $x_\varphi\in X$ such that $\<y,x_\varphi\> = \<\gamma_{y,x}, \varphi\>$ for all $y\in Y$.

For each $\varphi\in\mathcal{U}^{**}$, the map on $X$ that sends $x\mapsto x_\varphi$ is itself bounded and linear; thus we may define
$\rho(\varphi)x := x_{\varphi}$, and see that $\rho(\varphi)\in L(X)$. As $\norm{\rho(\varphi)}\leq \norm{\rho_0}\norm{\varphi}$ for 
all $\varphi\in\U^{**}$ and $\rho(f) = \rho_0(f)$ for all $f\in\U\subseteq\U^{**}$, it is clear that the norm of the map $\rho: \mathcal{U}^{**}\to L(X)$ is $\norm{\rho} = \norm{\rho_0}$.

By definition, $\<y,\rho(\varphi)x\> = \<\gamma_{y,x},\varphi\>$ for all $y\in Y, x\in X$, so we have that $\rho$ is linear. To see that $\rho$ is multiplicative we first verify its weak-$\ast$ continuity. Suppose $\varphi\in\U^{**}$ and take a net 
$\{\varphi_\alpha\}\subseteq\U^{**}$ such that for all $\gamma\in\U^{*}$,
$\<\gamma,\varphi\> = \lim_\alpha \< \gamma,\varphi_\alpha\>$.
Then for all $y\in Y, x\in X$, 
\begin{equation}\label{rhocont}
 \<y, \rho(\varphi)x\> = \< \gamma_{y,x}, \varphi\> 
	= \lim_\alpha \< \gamma_{y,x} ,\varphi_\alpha\> 
	= \lim_\alpha \<y, \rho(\varphi_\alpha)x\> 
\end{equation}
as desired. 

Now let $\varphi, \psi \in \mathcal{U}^{**}$. By Goldstine's Theorem, 
there exist nets $\{f_\alpha\}, \{g_\beta\} \subseteq \mathcal{U}$ such that 
$\<\gamma,\varphi\> = \lim_\alpha \< \gamma,f_\alpha\>\,,$ $\<\gamma,\psi\> = \lim_\beta \< \gamma,g_\beta\>$ for all $\gamma\in\U^*$.
We have that $\<\gamma,f_\alpha\star\psi\> = \lim_\beta \< \gamma,f_\alpha \star g_\beta\> = \lim_\beta \< \gamma,f_\alpha g_\beta\>$; 
and applying the continuity (\ref{rhocont}) gives
\[ \<y, \rho(f_\alpha\star\psi)x\> 	= \lim_\beta \< y,\rho(f_\alpha g_\beta)x\> 
%\\				&= \lim_\beta \< x,\rho_0(f_\alpha g_\beta)y\> \\
= \lim_\beta \< y,\rho_0(f_\alpha)\rho_0(g_\beta)x\> = \< y,\rho(f_\alpha)\rho(\psi)x\> \]
for all $y\in Y, x\in X$.
Similarly, we apply (\ref{rhocont}) to 
$\<\gamma,\varphi\star\psi\> = \lim_\alpha \< \gamma,f_\alpha\star \psi\>$, and use the above to see that 
$\<y, \rho(\varphi\star\psi)x\> = \< y,\rho(\varphi)\rho(\psi)x\>$ for all $y\in Y, x\in X$. Hence $\rho$ is multiplicative and we
have the desired homomorphism. 
%The uniqueness of $\rho$ is similarly derived from the continuity (\ref{rhocont}).
\end{proof}

We apply the theorem to the case at hand by noting that $\B(\sigma)$ is a subalgebra of $C(\sigma)^{**}$.

\begin{corollary}\label{contextcor}
Let $X$ be a Banach space such that $X = Y^*$ for some Banach space $Y$, and let $\sigma$ be a non-empty compact subset of $\C$. 
Then $T\in L(X)$ has a $\B(\sigma)$ functional calculus if and only if it has
a $C(\sigma)$ functional calculus.
\end{corollary}
\begin{proof}
If the $C(\sigma)$ functional calculus for $T$ is given by 
$\rho_0: C(\sigma)\rightarrow L(X)$, then from Theorem \ref{sillsextn}, $T$ has a $C(\sigma)^{**}$ functional calculus given by 
$\rho: C(\sigma)^{**} \to L(X)$. For each $\varphi\in \B(\sigma)$, the map $\mu \mapsto \int_\sigma \varphi\, d\mu$ ($\mu \in C(\sigma)^*$) 
defines a bounded linear functional on $C(\sigma)^*$. Thus we may consider $\varphi$ as an element of $C(\sigma)^{**}$ and write 
$\<\mu,\varphi\> = \int_\sigma \varphi\, d\mu$ for all $\mu\in C(\sigma)^*$, $\varphi\in \B(\sigma)$. A routine calculation verifies that 
the Arens multiplication in $C(\sigma)^{**}$ restricted to $\B(\sigma)$ coincides with the usual pointwise multiplication. Thus 
$\B(\sigma)$ with pointwise multiplication is a subalgebra of $C(\sigma)^{**}$ and 
the map $\trho : \B(\sigma) \to L(X)$, defined by $\trho = \rho|_{\B(\sigma)}$ gives the required functional calculus. 
\end{proof}

In \cite{nagy}, Nagy showed that possession of a $\B(\sigma)$ functional calculus is not only necessary, but also sufficient for an operator to be finitely scalar spectral. We know then, that an operator on a dual space with a $C(\sigma)$ functional calculus is finitely scalar spectral. However, we have more: in fact such an operator must be scalar-type prespectral.

\begin{proposition}\label{cor_prespec}
Let $X$ be a Banach space such that $X = Y^*$ for some Banach space $Y$.
Then $T\in L(X)$ is scalar-type prespectral of class $Y$ if and only if it has
a $C(\sigma)$ functional calculus for some non-empty compact set $\sigma\subseteq\C$.
\end{proposition}
\begin{proof}
If $T$ is scalar-type prespectral, then it has a $\B(\sigma(T))$ functional calculus \cite{dunfordspec} and thus a $C(\sigma(T))$ functional calculus. Conversely, if $T$ has a $C(\sigma)$ functional calculus, then by Corollary \ref{contextcor} we may construct a bounded homomorphism
$\trho : \B(\sigma) \to L(X)$ that defines a $\B(\sigma)$ functional calculus for $T$. Using $\borel$ to denote the Borel subsets of the plane, define $E:\borel\rightarrow L(X)$ by
$E(\Delta) = \trho(\chi_{\Delta\cap\sigma})$ for $\Delta\in\borel$.
Then $E(\sigma) = I$ and (directly from the properties of the characteristic function and the homomorphism)
$E$ is a finitely additive spectral measure. To show that $E$ has a weak countable additivity, take a countable collection of disjoint sets $\{\Delta_k\}\subseteq\borel$. We may assume that each $\Delta_k \subseteq \sigma$, as $E(\Delta) = 0$ whenever $\Delta\subseteq \C\take\sigma$. 
Let $\bigcup_{k=1}^{\infty}\Delta_k = \Delta$ and let $\nabla_n = \bigcup_{k=n+1}^{\infty} \Delta_k$. Then for $y\in Y, x\in X , n\in \N$,
\begin{eqnarray*}
 \<y,(E(\Delta) - \sum_{k=1}^{n}E(\Delta_{k}))x \>
		&=&  \<y,E(\nabla_n)x \> \\
		&=&  \<y,\trho(\chi_{\nabla_n})x \> \\
		&=&  \< \gamma_{y,x}, \chi_{\nabla_n} \> \\
		&=&  \int_{\sigma} \chi_{\nabla_n}\,d\gamma_{y,x} \,,
\end{eqnarray*}
where $\gamma_{y,x}\in C(\sigma)^*$ is defined as in the proof of Theorem \ref{sillsextn}. Letting $n\to\infty$ sends the integral to zero and we have that $E$ is a spectral measure of class $Y$.

To verify that $E$ is a resolution of the identity for $T$, firstly observe that as
$T = \trho(\lam)$, $T$ commutes with $E(\Delta) = \trho(\chi_\Delta)$ 
for all $\Delta\in\borel$. Basic techniques of local spectral theory now show that for any $\Delta\in\borel$, where $X_T(\Delta)$ denotes the local spectral subspace for $T$, we have
\[ E(\Delta)X = \trho(\chi_\Delta)X \subseteq X_T(\supp(\chi_\Delta)) = X_T(\overline{\Delta})\,. \]
Thus $\sigma(T|E(\Delta)X) \subseteq \overline{\Delta}$.
\end{proof}

This result provides a companion to a well known result of Berkson and Dowson \cite{berkdowprespec} stating that if an operator $T\in L(Y)$ has a continuous functional calculus, then $T^*\in L(Y^*)$ is scalar-type prespectral of class $Y$. A corollary of that result is that the adjoint of a finitely spectral operator is a prespectral operator. Proposition \ref{cor_prespec} also extends easily to non-scalar operators.

\begin{corollary}\label{coincide}
On a dual Banach space the classes of finitely spectral and prespectral operators coincide. 
\end{corollary}

We note that it is still an open question as to whether there exists an operator with a $\B(\sigma)$ functional calculus (that is, a finitely spectral operator) that is not prespectral of any class. It is known (see \cite{nagy} and \cite{gillespie}) that on a Banach space not
containing a subspace isomorphic to $\ell^\infty$, the classes of finitely spectral, prespectral and spectral operators 
coincide, so combining this with Corollary \ref{coincide} leaves us with a rather restricted choice of spaces on which to 
search for such an operator.

\section{Real $\AC(\sigma)$ functional calculus}\label{section_acextns}

We now look at some applications of the above results to real $\AC(\sigma)$ operators, that is, operators admitting an $\AC(\sigma)$ functional calculus for a non-empty compact set $\sigma\subseteq\R$. 
Real $\AC(\sigma)$ operators have traditionally been referred to as well-bounded operators, and we refer the reader to \cite{dowson} and \cite{ashdoust1} for details of their theory. As in the case of operators with a $C(\sigma)$ functional calculus, much is known about real $\AC(\sigma)$ operators acting on reflexive spaces. A real $\AC(\sigma)$ functional calculus for an operator $T$ on a reflexive space $X$ is necessarily weakly compact. Unlike a weakly compact $C(\sigma)$ functional calculus, this does not provide a spectral measure for $T$. However it does ensure the existence of a `spectral family' of projections $\{E(\lambda)\}_{\lambda\in\R}$ that act on the space $X$ and provide an integral representation for $T$ \cite{berkdow}. The integral representation can be used to extend the $\AC(\sigma)$ functional calculus to $\BV(\sigma)$, the algebra of functions of bounded variation defined on $\sigma$.

If $X$ is a non-reflexive Banach space, a real $\AC(\sigma)$ operator $T\in L(X)$ is only guaranteed an associated (non-unique) family of projections $\{F(\lambda)\}_{\lambda \in \R}$ that act on the dual space $X^*$, and represent $T$ via
\begin{equation}\label{dotirep}
 \<Tx,x^*\> = b\<x,x^*\> - \int_a^b \<x,F(\lambda)x^*\> \, d\lambda,\; x \in X, x^* \in X^* 
\end{equation}
for some $[a,b]\subseteq\R$. The family $\{F(\lambda)\}$ is called a decomposition of the identity for $T$; if the family is unique, $T$ is said to be uniquely decomposable. 
Now the representation (\ref{dotirep}) may not, in general, be used to construct a $BV(\sigma)$ functional calculus for an operator. However, many real $\AC(\sigma)$ operators on non-reflexive spaces do admit a $\BV(\sigma)$ functional calculus. Thus a natural question is to ask for conditions on a space $X$ or an operator $T$ that ensure such an extended functional calculus exists. In \cite{doustT}, it is shown that if a space $X$ contains a complemented copy of either $c_0$ or $\ell^1$, then there exists a real $\AC(\sigma)$ operator $T\in L(X)$ that does not admit a $\BV(\sigma)$ functional calculus. As $\ell^1$ is a dual space, this means that we cannot hope for a direct analogue of Corollary \ref{contextcor} here.  

In fact, it is not known if there exists a non-reflexive Banach space $X$ that can guarantee a $\BV(\sigma)$ functional calculus for all 
real $\AC(\sigma)$ operators $T\in L(X)$. 
Following Sills \cite{sills}, we may apply Theorem \ref{sillsextn} and the fact that the algebra $\BV[a,b]$ is isomorphic to a quotient algebra of $\AC[a,b]^{**}$ to find that a real $\AC(\sigma)$ operator $T$ on a dual space $X$ has a $\BV(\sigma)$ functional calculus under certain conditions. One sufficient condition is that $T$ possess a decomposition of the identity $\{F(\lambda)\}$ of bounded variation \cite{teraudsthesis}, that is, that the function $\lambda\mapsto F(\lambda)$ be of bounded variation on $\R$. We omit the details of this here, as a stronger result can easily be gained by exploiting the relationship between continuous and absolutely continuous functional calculi.

The algebra $\AC(\sigma)$ is equipped with the norm $\norm{\cdot}_{\BV}$ (see \cite{ashdoust1}), and for every 
$f\in \AC(\sigma)$, we have $\norm{f}_\infty \leq \norm{f}_{\BV}$. Thus every operator with a $C(\sigma)$ functional calculus 
is an $\AC(\sigma)$ operator, and every operator with a $\B(\sigma)$ functional calculus admits a $\BV(\sigma)$ one. 
This is true for any compact $\sigma\subseteq\C$, but when the set $\sigma$ is real, Berkson and Dowson \cite[Theorem 5.2]{berkdow} showed that  for an operator $T \in L(X)$, possession of a $C(\sigma)$ functional calculus is both necessary and sufficient for $T$ to be a real $\AC(\sigma)$ operator with a decomposition of the identity of bounded variation. Combining this with Corollary \ref{contextcor} and Proposition \ref{cor_prespec}, we have the following. 

\begin{theorem}
Let $X$ be a Banach space such that $X = Y^*$ for some Banach space $Y$ and $\sigma$ be a non-empty compact subset of $\R$. 
Then for $T\in L(X)$, the following are equivalent.
\begin{enumerate}
\item[(i)] $T$ is a real $\AC(\sigma)$ operator with a decomposition of the identity of bounded variation
\item[(ii)] $T$ has a $C(\sigma)$ functional calculus
\item[(iii)] $T$ has a $\B(\sigma)$ functional calculus
\item[(iv)] $T$ is scalar-type prespectral of class $Y$
\end{enumerate}
If any of the above hold, then the real $\AC(\sigma)$ operator $T$ is uniquely decomposable.
\end{theorem}

We note that a very similar result holds on any Banach space containing no copy of $c_0$ (see \cite{berkdow} and \cite{DdeLaub}), although in such a case an operator satisfying any of the conditions (i) - (iii) will be scalar-type spectral.


\begin{thebibliography}{20cm}

\bibitem{ashdoust1}
B. Ashton, and I. Doust, Functions of bounded variation on compact
subsets of the plane, Studia Math., \textbf{169} (2005), 163--188.

\bibitem{doustT}
I. Doust and V. Terauds, Extensions of an $AC(\sigma)$ functional
calculus, submitted.

\bibitem{berkdow}
E. Berkson\ and\ H. R. Dowson, On uniquely decomposable well-bounded operators, 
{\em Proc. London Math. Soc.} (3) {\bf 22} (1971), 339--358.

\bibitem{berkdowprespec}
E. Berkson\ and\ H. R. Dowson, Prespectral operators, Illinois J. Math. {\bf 13} 
(1969), 291--315.

\bibitem{DdeLaub} 
I. Doust\ and\ R. deLaubenfels, Functional calculus, integral representations, and Banach 
space geometry, {\em Quaestiones Math}. {\bf 17} (1994), no.~2, 161--171.

\bibitem{dowson}
H. R. Dowson, {\it Spectral theory of linear operators}, Academic Press, London, 1978.

\bibitem{dunfordspec}
N. Dunford, Spectral operators, Pacific J. Math. {\bf 4} (1954), 321--354.

\bibitem{dunford3}
N. Dunford\ and\ J. T. Schwartz, {\it Linear operators. Part III}, Interscience Publishers 
[John Wiley \& Sons, Inc.], New York, 1971.

\bibitem{gillespie} 
T. A. Gillespie, Spectral measures on spaces not containing $l\sp{\infty }$, 
Proc. Edinburgh Math. Soc. (2) {\bf 24} (1981), no.~1, 41--45.

\bibitem{nagy}
B. Nagy, Finitely spectral operators, Glasgow Math. J. {\bf 28} (1986), no.~1, 95--112. 

\bibitem{sills}
W. H. Sills, On absolutely continuous functions and the well-bounded operator, 
{\em Pacific J. Math.} {\bf 17} (1966), 349--366. 

\bibitem{spain}
P. G. Spain, On scalar-type spectral operators, Proc. Cambridge Philos. Soc. {\bf 69} 
(1971), 409--410.

\bibitem{teraudsthesis}
V. Terauds, Extensions of functional calculus for Banach space operators, PhD thesis, 
University of New South Wales (2006).

\end{thebibliography}
\end{document}